\documentstyle[11pt]{article}
 \newcommand{\qb}[2]{{\left [{#1 \atop #2} \right]}}
 \newlength{\standardunitlength}
\newtheorem{cor}{Corollary} \newtheorem{lemma}{Lemma}
\newtheorem{theorem}{Theorem} \newtheorem{prop}{Proposition}
\newenvironment{proof}{\noindent {\sc Proof:}}{$\Box$ \vspace{2 ex}}

\begin{document}

\begin{center} Applications of the Brauer complex: card shuffling,
permutation statistics, and dynamical systems \end{center}

\begin{center}
By Jason Fulman
\end{center}

\begin{center}
Stanford University
\end{center}

\begin{center}
Department of Mathematics
\end{center}

\begin{center}
Building 380, MC 2125
\end{center}

\begin{center}
Stanford, CA 94305
\end{center}

\begin{center}
To appear in J. Algebra
\end{center}

\begin{center}
email:fulman@math.stanford.edu
\end{center}

1991 AMS Primary Subject Classifications: 20G40, 20F55

\newpage Proposed running head: Applications of the Brauer

\newpage \begin{abstract} By algebraic group theory, there is a map
from the semisimple conjugacy classes of a finite group of Lie type to
the conjugacy classes of the Weyl group. Picking a semisimple class
uniformly at random yields a probability measure on conjugacy classes
of the Weyl group. Using the Brauer complex, it is proved that this
measure agrees with a second measure on conjugacy classes of the Weyl
group induced by a construction of Cellini using the affine Weyl
group. Formulas for Cellini's measure in type $A$ are found. This
leads to new models of card shuffling and has interesting
combinatorial and number theoretic consequences. An analysis of type
$C$ gives another solution to a problem of Rogers in dynamical
systems: the enumeration of unimodal permutations by cycle
structure. The proof uses the factorization theory of palindromic
polynomials over finite fields. Contact is made with symmetric
function theory. \end{abstract}


\section{Introduction}

	In performing a definitive analysis of the
Gilbert-Shannon-Reeds model of card-shuffling, Bayer and Diaconis
\cite{BD} defined a one-parameter family of probability measures on
the symmetric group $S_n$ called $k$-shuffles. Given a deck of $n$
cards, one cuts it into $k$ piles with probability of pile sizes
$j_1,\cdots,j_k$ given by $\frac{{n \choose
j_1,\cdots,j_k}}{k^n}$. Then cards are dropped from the packets with
probability proportional to the pile size at a given time (thus if the
current pile sizes are $A_1,\cdots,A_k$, the next card is dropped from
pile $i$ with probability $\frac{A_i}{A_1+\cdots+A_k}$). They proved
that $\frac{3}{2}log_2(n)$ $2$-shuffles are necessary and suffice to
mix up a deck of $n$ cards. Aldous \cite{A} had previously obtained
this bound asymptotically in $n$; the paper \cite{F4} shows that the
use of cuts does not help to speed things up.

	One motivation for the current paper is the fact that GSR
measures are well-studied and appear in many mathematical
settings. The paper \cite{H} is a good reference for applications to
Hochschild homology (tracing back to Gerstenhaber and Schack
\cite{Ger}), and the paper \cite{BergeronWolfgang} describes the
relation with explicit versions of the Poincar\'e-Birkhoff-Witt
theorem. Section 3.8 of \cite{SSt} describes GSR shuffles in the
language of Hopf algebras. In recent work, Stanley \cite{Sta} has
related biased riffle shuffles with the Robison-Schensted-Knuth
correspondence, thereby giving an elementary probabilistic
interpretation of Schur functions and a different approach to some
work of interest to the random matrix community. He recasts many of
the results of \cite{BD} and \cite{F} using quasisymmetric
functions. Connections of riffle shuffling with dynamical systems
appear in \cite{BD}, \cite{La1}, \cite{La2}. Generalizations of the
GSR shuffles to other Coxeter groups, building on the papers
\cite{BBHT},\cite{BHR} appear in \cite{F2} and \cite{F3}.

	For further motivation, it is useful to recall one of the most
remarkable yet mysterious properties of these $k$-shuffles. Since
$k$-shuffles induce a probability measure on conjugacy classes of
$S_n$, they induce a probability measure on partitions $\lambda$ of
$n$. Consider the factorization of random degree $n$ polynomials over
a field $F_q$ into irreducibles. The degrees of the irreducible
factors of a randomly chosen degree $n$ polynomial also give a random
partition of $n$. The fundamental result of Diaconis-McGrath-Pitman
(DMP) \cite{DMP} is that this measure on partitions of $n$ agrees with
the measure induced by card shuffling when $k=q$. This allowed natural
questions on shuffling to be reduced to known results on factors of
polynomials and vice versa.

	There are three different proofs of the DMP result, each of
them mysterious in its own way. The first proof, in \cite{DMP}, is
combinatorial and makes use of a magical bijection of Gessel and
Reutenauer \cite{G}; they include a self-contained proof of this
bijection. The second proof, in \cite{H}, proves an equivalent
assertion about induced characters, but also uses the
Gessel/Reutenauer bijection. (The equivalence between the DMP theorem
and the results on induced characters is not completely obvious; see
Section 4 of \cite{F4} for explanation). The third, and perhaps most
principled proof of the equivalent assertion about induced characters
appears in \cite{BBGar}, using facts about free Lie algebras
\cite{Gar}. But it is unclear how to generalize free Lie Algebras to
arbitrary types. Motivated by the observation that degree $n$
polynomials over $F_q$ are the semisimple orbits of $GL(n,q)$ on its
Lie algebra, the paper \cite{F3} gave Lie theoretic reformulations and
generalizations of the DMP result. However the proofs still strike one
as unnatural, and it is not clear when the generalizations hold. We
should in fairness point out that there is hope of a uniform
generalization of the Gessel/Reutenauer bijection, at least for the
conjugacy class of Coxeter elements \cite{Rei1},\cite{Rei2} and that
there is an analog of the Free Lie algebra in type $B$ \cite{B}.

	The goal of the current article is to study a setting where
these complications vanish and an analog (Theorem \ref{analogDMP}) of
the DMP result holds in all types, and can be proved in a uniform and
natural way. The key idea is to study semisimple conjugacy classes in
groups such as $SL(n,q)$ rather than $GL(n,q)$; then these polynomials
can be viewed as points in Euclidean space, and the extra geometric
structure forces natural choices.

	The precise contents of this paper are as follows. Section
\ref{background} begins by describing the algebraic groups set-up and
the map $\Phi$ from semisimple conjugacy classes to conjugacy classes
of the Weyl group, giving examples. It makes a connection with the
Gessel/Reutenauer map, demystifying it somewhat. Section
\ref{background2} gives a probabilistic version of a construction of
Cellini \cite{Ce1} and states the analog of the DMP theorem, which
will be proved in Section \ref{refine}.

	Section \ref{typeA} focuses on understanding Cellini's
construction in types $A$. It emerges that in type $A$ the probability
of a permutation involves both its number of cyclic descents and major
index. This is interesting because while combinatorialists have
thoroughly studied the joint distribution of permutations by descents,
major index, and cycle structure \cite{Ge}, problems involving cyclic
descents have not been treated and are regarded by the experts as
harder. It is also shown that the type $A$ construction leads to new
models of card shuffling. Section \ref{typeA} then shows that even for
the identity conjugacy class in type $A$, Theorem \ref{analogDMP}
gives an interesting result--a number theoretic reciprocity law. For
more general conjugacy classes, Theorem \ref{analogDMP} is given a
formulation in terms of generating functions which highlights the
connections with number theory.

	Section \ref{typeC} studies Cellini's construction in type
$C$. Unlike the type $A$ case, formulas for type $C$ follow easily
from work of Cellini (this is one of the few times in mathematics when
the hyperoctahedral group is easier to understand than the symmetric
group). Thus the main point of this section is to give interpretations
in terms of card shuffling. This unifies work of \cite{BD} and
\cite{BB} and implies a simple formula for Bayer-Diaconis
hyperoctahedral shuffles. Some work of Bob Beals on total variation
distance of hyperoctahedral shuffles to uniform is understood is a new
way.

	Section \ref{dynamical} gives applications of the type $C$
analog of the DMP theorem to dynamical systems. Specifically, an
alternate solution to a problem posed by Rogers \cite{Ro}--the
enumeration of unimodal permutations by cycle structure--is given. The
mathematics in this section was accepted by J. Algebra in 1/01, prior
to the appearance of the preprint \cite{T}; hence Section
\ref{dynamical} gives the first derivation of the cycle index of
unimodal permutations. (The proof in \cite{T}, however, gives the
first derivation using only symmetric functions and the
Gessel-Reutenauer machine, generalizing reasoning here for
$n$-cycles). The cycle index leads to interesting asymptotic
results. The section closes by giving a more conceptual proof of a
combinatorial result of Reiner. Further generalizations of the DMP
theorem based on symmetric functions are in \cite{F5}.

	Section \ref{refine} uses the Brauer complex (an object
originally introduced in modular representation theory) to prove
Theorem \ref{analogDMP}, conjectured in an early version of this
paper. The proof presented here is definitive and was provided by
Professor Roger Carter in October 2000; he selflessly declined to be a
coauthor. His work supersedes the author's type $C$ odd characteristic
proof and unsuccessful efforts to exploit the viewpoint of semisimple
conjugacy classes as points in Euclidean space; these partial results
remain at http://xxx.lanl.gov/ in the paper ``Cellini's ...'' but have
been cut from this final version.

\section{Algebraic Groups} \label{background}

	Notation about algebraic groups will conform to that in
$\cite{C}$, which together with \cite{Hu1} contains all the relevant
background for this paper. The book \cite{Hu2} is a good reference for
information about Coxeter groups. Throughout $G$ is a semisimple,
simply connected group defined over a finite field of $q$
elements. Letting $F$ be a Frobenius automorphism of $G$, we suppose
that $G$ is $F$-split. Pages 39-41 of \cite{C} list these groups. For
instance in type $A_{n-1}$ the group is $SL(n,q)$ and in type $C_n$
the group is $Sp(2n,q)$. Let $W$ denote the Weyl group. This is the
symmetric group in the first example and the hyperoctahedral group in
the second example.

	There is a natural map $\Phi$ from semisimple conjugacy
classes $c$ of $G^F$ to conjugacy classes of the Weyl group. Let $x$
be an element in the class $c$. Theorem 2.11 of \cite{Hu1} implies
that the centralizers of semisimple elements of $G$ are
connected. Consequently $C_G(x)$, the centralizer in $G$ of $x$, is
determined up to $G^F$ conjugacy. As is possible from page 33 of
\cite{C}, let $T$ be an $F$-stable maximally split maximal torus in
$C_G(x)$; $T$ is determined up to $G^F$ conjugacy. Proposition 3.3.3
of \cite{C} gives that the $G^F$ conjugacy classes of $F$-stable
maximal tori of $G$ are in bijection with conjugacy classes of
$W$. Define $\Phi(c)$ to be the corresponding conjugacy class of $W$.

	Lemma \ref{explicit} makes the map $\Phi$ explicit. Recall
from Proposition 3.7.3 of \cite{C} that there is a bijection between
the semisimple conjugacy classes of $G^F$ and the $F$-stable orbits in
$T/W$.

\begin{lemma} \label{explicit} Let $T_0$ be a maximally split torus of
$G$ and let $t_0 \in T_0$ be a representative for the semisimple
conjugacy class $c$. Let $\Psi_0$ be the root system of $C_G(t_0)$,
i.e. all roots $\alpha$ such that $\alpha(t_0)=1$. Suppose that
$F(t_0)=t_0^w$ and that $w^{-1}(\Psi_0^+)=\Psi_0^+$ for some positive
system of $\Psi_0$. Then $\Phi(c)$ is the conjugacy class of
$w$. \end{lemma}

\begin{proof} Let $T$ be an $F$-stable maximal torus of $G$ obtained
by twisting $T_0$ by $w$. Let $t$ be the image of $t_0$ under the
corresponding conjugation map. Then $T$ is a maximally split torus in
$C(t)$ if and only if there is an $F$-stable Borel subgroup of $C(t)$
containing $T$, which happens if and only if there is a positive
system of $\Psi_0$ such that $w^{-1}(\Psi_0^+)=\Psi_0^+$. \end{proof}

	Next we give two examples, which will be used later in this
paper.

\begin{enumerate}

\item The first example is $SL(n,q)$ with Weyl group the symmetric
group on $n$ symbols. The semisimple conjugacy classes $c$ correspond
to monic degree $n$ polynomials over $F_q$ with constant term
$1$. Such a polynomial factors into irreducible polyomials. Let $n_i$
be the number (counted with multiplicity) of these irreducible factors
of degree $i$. Then the corresponding conjugacy class $\Phi(c)$ has
$n_i$ cycles of length $i$.

	To see this from the lemma, note that $t$ takes the form of a
diagonal matrix, where the entries along the diagonal are the roots of
the characteristic polynomial in the algebraic closure. The Frobenius
map $F$ acts by raising elements to the $q$th power, thus permuting
the elements along the diagonal. (This permutation is unique if for
instance all of the irreducible factors of the characteristic
polynomial are distinct. If the irreducible factors aren't all
distinct, then the root system $\Psi_0$ is non-trivial, since the root
$e_i-e_j$ sends a diagonal matrix with diagonal entries
$(x_1,\cdots,x_n)$ to $x_i/x_j$, and hence some roots send $t$ to
$1$). If one considers the particular positive subsystem $\Psi_0 \cap
\{e_i-e_j: i <j\}$ of $\Psi_0$, it is easy to see explicitly that
there is a unique $w$ satisfying $F(t_0)=t_0^w$,
$w^{-1}(\Psi_0^+)=\Psi_0^+$ and that $w$ has $n_i$ cycles of length
$i$.

	We remark that the map $\Phi$ is closely related to the
Gessel/Reutenauer map. The Gessel/Reutenauer map associates to each
multiset of primitive (i.e. not equal to any of its proper rotations)
necklaces on the symbols $\{0,1,\cdots,q-1\}$ a permutation $w$. This
map is carefully exposited in \cite{G} and was used in shuffling work
in \cite{DMP}; we omit its definition here.

	If one fixes generators for the multiplicative group of each
finite extension of $F_q$, the monic degree $n$ polynomials $\phi$
correspond to multisets of primitive necklace \cite{Go}. For example
suppose that $\theta$ is a generator of the multiplicative group of
$F_{q^3}$. Then a degree $3$ monic irreducible polynomial corresponds
to an orbit of $\theta^i$ under the Frobenius map, for some
$i$. Writing this $i$ base $q$ gives a size 3 primitive necklace. To
each necklace entry one can associate a root of $\phi$, by taking
$\theta^{j}$ where $j$ is a number base $q$ obtained by rotating $i$
so that the specified necklace entry is the leftmost digit of $j$. For
example if the necklace is $(01011)$ and one is working base $2$ then
the middle $0$ would correspond to $\theta^j$ where $j$ is $01101$
base $2$. For what follows it is helpful to make $01101$ an infinite
word by repeating it $0110101101...$ Now associated to $\phi$ one can
form a diagonal matrix whose elements are the roots of $\phi$, ordered
lexicographically by their associated infinite word. Then the
permutation associated to this matrix through Lemma \ref{explicit} is
equal to the permutation which the Gessel-Reutenauer map associated to
the corresponding multiset of primitive necklaces.

\item The second example is type $C$. The group in question is
$Sp(2n,q)$ with Weyl group $C_n$ the group of signed permutations. The
semisimple conjugacy classes $c$ correspond to monic degree $2n$
polynomials $\phi(z)$ with non-zero constant term that are invariant
under the involution sending $\phi(z)$ to $\bar{\phi}(z) =
\frac{z^{2n} \phi(1/z)}{\phi(0)}$. Such polynomials can be more simply
described as monic degree $2n$ polynomials which are palindromic in
the sense that the coefficient of $z^i$ is equal to the coefficient of
$z^{2n-i}$. These factor uniquely into irreducibles as

\[ \prod_{\{\phi_j,\bar{\phi_j}\}} [\phi_j \bar{\phi_j}]^{r_{\phi_j}}
\prod_{\phi_j: \phi_j=\bar{\phi_j}} \phi_j^{s_{\phi_j}}\] where the
$\phi_j$ are monic irreducible polynomials and $s_{\phi_j} \in
\{0,1\}$. The conjugacy classes of $C_n$ correspond to pairs of
vectors $(\vec{\lambda},\vec{\mu})$ where
$\vec{\lambda}=(\lambda_1,\cdots,\lambda_n)$,
$\vec{\mu}=(\mu_1,\cdots,\mu_n)$ and $\lambda_i$ (resp. $\mu_i$) is
the number of positive (resp. negative) $i$ cycles of an element of
$C_n$, viewed as a signed permutation. From Lemma \ref{explicit} one
can see that the conjugacy class of $C_n$ corresponding to $c$ is then
determined by setting $\lambda_i=\sum_{\phi: deg(\phi)=i} r_{\phi}$
and $\mu_i=\sum_{\phi: deg(\phi)=2i} s_{\phi}$.

\end{enumerate}

\section{Cellini's Work} \label{background2}

	Next we recall the work of Cellini \cite{C} (the definition
which follows differs slightly from hers, being inverse, making use of
her Corollary 2.1, and renormalizing so as to have a probability
measure). We follow her in supposing that $W$ is a Weyl group (i.e. a
finite reflection group which arises from a Chevalley group). Let
$\Delta=\{\alpha_1,\cdots,\alpha_r\}$ be a simple root system for
$W$. Letting $\alpha_0$ denote the negative of the highest root, let
$\tilde{\Pi}=\Pi \cup \alpha_0$. Define the cyclic descent $Cdes(w)$
to be the elements of $\tilde{\Pi}$ mapped to negative roots by $w$,
and let $cd(w)=|Cdes(w)|$. For future use we remark that the descent
set of $w$ is defined as the subset of $\Delta$ mapped to negative
roots by $w$.

	For instance for $S_n$ the simple roots with respect to a
basis $e_1,\cdots,e_n$ are $e_i-e_{i+1}$ for $i=1,\cdots,n-1$ and
$\alpha_0=e_n-e_1$. Thus the permutation $4 \ 1 \ 3 \ 2 \ 5$ (in
2-line form) has 3 cyclic descents and 2 descents. Type $C$ examples
will be treated in Section \ref{typeC}.

	Now we use cyclic descents to define shuffles. For $I
\subseteq \tilde{\Pi}$, put \[ U_I = \{w \in W|Cdes(w) \cap I =
\emptyset\} .\] Let $Y$ be the coroot lattice. Then define $a_{k,I}$
by

\[ |\{ y \in Y| <-\alpha_0,y> = k, <\alpha_i,y> = 0 \ for \ \alpha_i
\in I-\alpha_0, <\alpha_i,y> >0 \ for \ \alpha_i \in \tilde{\Pi}-I\}|
\ if \ \alpha_0 \in I \]

\[ |\{ y \in Y| <-\alpha_0,y> < k, <\alpha_i,y> = 0 \ for \ \alpha_i
\in I, <\alpha_i,y> >0 \ for \ \alpha_i \in \Pi-I\}| \ if
 \ \alpha_0 \not \in I.\] Finally, define an element $x_k$ of the
group algebra of $W$ by \[ x_k = \frac{1}{k^r} \sum_{I \subseteq
\tilde{\Pi}} a_{k,I} \sum_{w \in U_I} w .\] Equivalently, the
coefficient of an element $w$ in $x_k$ is \[ \frac{1}{k^r} \sum_{I
\subseteq \tilde{\Pi} - Cdes(w)} a_{k,I}.\] This coefficient will be
denoted by $x_k(w)$ throughout the paper. We will refer to these $x_k$
as Affine $k$-shuffles.

	In type $A_{n-1}$ this says that the coefficient of $w$ is
$x_k$ is equal to $\frac{1}{k^{n-1}}$ multiplied by the number of
integers vectors $(v_1,\cdots,v_n)$ satisfying the conditions

\begin{enumerate}
\item $v_1+\cdots+v_n=0$
\item $v_1 \geq v_2 \geq \cdots \geq v_n, v_1-v_n \leq k$
\item $v_i>v_{i+1}$ if $w(i)>w(i+1)$ (with $1 \leq i \leq n-1$)
\item $v_1<v_n+k$ if $w(n)>w(1)$
\end{enumerate}

	From Cellini (loc. cit.), it follows that the $x_k$ satisfy
the following two desirable properties:

\begin{enumerate}
\item (Measure) The sum of the coefficients in the expansion of $x_k$
in the basis of group elements is 1. Equivalently,

\[ \sum_{I \subseteq \tilde{\Pi}} a_{k,I} |U_I| = k^r .\] In
probabilistic terms, the element $x_k$ defines a probability measure
on the group $W$.

\item (Convolution) $x_k x_h = x_{kh}$.
\end{enumerate}

	The above definition of $x_k$ is computationally convenient
for this paper. We note that Cellini (loc. cit.) constructed the $x_k$
in the following more conceptual way, when $k$ is a positive
integer. Let $W_k$ be the index $k^r$ subgroup of the affine Weyl
group that is generated by reflections in the hyperplanes
corresponding to $\{\alpha_1,\cdots,\alpha_r\}$ and also the
hyperplane $\{<x,-\alpha_0>=k\}$. There are $k^r$ unique minimal
length coset representatives for $W_k$ in the affine Weyl group, and
$x_k$ is obtained by projecting them to the Weyl group.

	The following problem is very natural. We remark that for GSR
riffle shuffles, Problem 1 was studied by Hanlon \cite{H}. Diaconis
has been a vigorous advocate of such questions, emphasizing the link
with convergence rates of Markov chains \cite{D}.

{\it Problem 1:} Determine the eigenvalues (and multiplicities) of
$x_k$ acting on the group algebra by left multiplication. More
generally, recall that the Fourier transform of a probability measure
$P$ at an irreducible representation $\rho$ is defined as $\sum_{w \in
W} P(w) \rho(w)$. For each $\rho$, what are the eigenvalues of this
matrix?

	To close the section, we state the analog of the DMP theorem.

\begin{theorem} \label{analogDMP} Let $G$ be a semisimple, simply
connected group defined over a finite field of $q$ elements. Letting
$F$ be a Frobenius automorphism of $G$, suppose that $G$ is
$F$-split. Let $c$ be a semisimple conjugacy class of $G^F$ chosen
uniformly at random. Then for all conjugacy classes $C$ of the Weyl
group $W$, \[ \sum_{w \in C} Probability(\Phi(c) = C) = \sum_{w \in C}
Coef. \ of \ w \ in \ x_q. \] \end{theorem}

\section{Type $A$ Affine Shuffles} \label{typeA}

	To begin, we derive four expressions for $x_k$ in type
$A_{n-1}$. For this recall that the major index of $w$ is defined by
$maj(w) = \sum_{i: 1 \leq i \leq n-1 \atop w(i)>w(i+1)} i$. It is the
sum of the positions of the descents of $w$. The notation $\qb{n}{k}$
denotes the $q$-binomial coefficient
$\frac{(q^n-1)\cdots(q-1)}{(q^k-1)\cdots (q-1) (q^{n-k}-1) \cdots
(q-1)}$. Let $C_m(n)$ denote the Ramanujan sum $\sum_k e^{\frac{2\pi
ikn}{m}}$ where $k$ runs over integers prime to $m$ satisfying $1 \leq
k \leq m$.

	The following lemma of Von Sterneck (see \cite{Ram} for a
proof in English) will be helpful. We emphasize that it is only used
in the derivation of the fourth formula for $x_k$; the first three
expressions will not require it.

\begin{lemma} \label{VS} (\cite{V}) The number of ways of expressing
$n$ as the sum mod $m$ of $k \geq 1$ integers of the set
$0,1,2,\cdots,m-1$ repetitions being allowed is \[ \frac{1}{m}
\sum_{d|m,k} {\frac{m+k-d}{d} \choose \frac{k}{d}} C_d(n).\]
\end{lemma}

	Recall that $x_k(w)$ is the coefficient of $w$ in $x_k$. For
us partitions have the standard number theoretic meaning as in
\cite{HW}.

\begin{theorem} \label{typeAform1} In type $A_{n-1}$, $x_k(w)$ is
equal to any of the following:

\begin{enumerate}

\item $\frac{1}{k^{n-1}}$ multiplied by the number of partitions with
$\leq n-1$ parts of size at most $k-cd(w)$ such that the total number
being partitioned has size congruent to $-maj(w) \ mod \ n$.

\item $\frac{1}{k^{n-1}}$ multiplied by the number of partitions with
$\leq k-cd(w)$ parts of size at most $n-1$ such that the total number
being partitioned has size congruent to $-maj(w) \ mod \ n$.

\item \[ \frac{1}{k^{n-1}} \sum_{r=0}^{\infty} Coeff. \ of \ q^{r \cdot n} \
in \ \left(q^{maj(w)} \qb{k+n-cd(w)-1}{n-1}\right).\]

\item \[ \begin{array}{ll} \frac{1}{n k^{n-1}} \sum_{d|n,k-cd(w)} {\frac{n+k-cd(w)-d}{d} \choose
\frac{k-cd(w)}{d}} C_d(-maj(w)) & \mbox{if \ $k-cd(w)>0$}\\
\frac{1}{k^{n-1}}  & \mbox{if $k-cd(w)=0, maj(w)=0$ mod n}\\
0 & \mbox{otherwise}
\end{array} \]

\end{enumerate} \end{theorem}

\begin{proof} From the definition of $x_k$,

\begin{eqnarray*} x_k(w) & = & \frac{1}{k^{n-1}} \sum_{I \subseteq
\tilde{\Pi}-Cdes(w)} a_{k,I}\\ & = & \frac{1}{k^{n-1}}
\sum_{v_1+\cdots+v_n=0, v_1 \geq \cdots v_{n-1} \geq v_n, v_1-v_n \leq
k, \vec{v} \in Z^n, \atop v_i > v_{i+1} \ if \ e_i-e_{i+1} \in
Cdes(w), \ and \ v_1-v_n<k \ if \ \alpha_0 \in Cdes(w)} 1\\ & = &
\frac{1}{k^{n-1}} Coeff. \ of \ q^{0} \ in \sum_{v_1 \geq \cdots
v_{n-1} \geq v_n, v_1-v_n \leq k, \vec{v} \in Z^n, \atop v_i > v_{i+1}
\ if \ e_i-e_{i+1} \in Cdes(w), \ and \ v_1-v_n<k \ if \ \alpha_0 \in
Cdes(w)} q^{v_1+\cdots+v_n}.  \end{eqnarray*} Given a vector
$(v_1,\cdots,v_n)$ with $\sum_i v_i=0$, one can translate it to
$(v_1-v_n,v_2-v_n,\cdots,0)$. The last coordinate of the new vector is
$0$ and the sum of the coordinates in this new vector is a multiple of
$n$. Abusing notation, we call this new vector
$(v_1,\cdots,v_n)$. Thus

\begin{eqnarray*}
x_k(w) & = & \frac{1}{k^{n-1}}
\sum_{r=0}^{\infty} Coeff. \ of \ q^{r \cdot n} \ in \sum_{ k \geq v_1
\geq \cdots v_{n-1} \geq v_n = 0, \vec{v} \in Z^n \atop v_i > v_{i+1}
\ if \ e_i-e_{i+1} \in Cdes(w), \ and \ v_1<k \ if \ 
\alpha_0 \in Cdes(w)} q^{\sum v_i}.
\end{eqnarray*} Now let $v_i'=v_i - |\{j: i \leq j \leq n-1, w(j)>w(j+1)\}|$. Then
the expression for $x_k(w)$ simplifies to

\begin{eqnarray*} 
&& \frac{1}{k^{n-1}}
\sum_{r=0}^{\infty} Coeff. \ of \ q^{r \cdot n} \ in \sum_{ k-cd(w) \geq v_1'
\geq \cdots v_{n-1}' \geq v_n' = 0, \vec{v} \in Z^n} q^{\sum v_i' + \sum_i |\{j: i \leq j \leq n-1, w(j)>w(j+1)\}|}\\
& = &  \frac{1}{k^{n-1}}
\sum_{r=0}^{\infty} Coeff. \ of \ q^{r \cdot n} \ in \ q^{maj(w)} \sum_{ k-cd(w) \geq v_1'
\geq \cdots v_{n-1}' \geq v_n' = 0, \vec{v} \in Z^n} q^{\sum v_i'}.
\end{eqnarray*}

	This proves the first assertion of the theorem. The second
assertion follows from the first by viewing partitions diagramatically
and taking transposes. The third assertion follows from either the
first or second assertions together with the well-known fact that the
generating function for partitions with at most $a$ parts of size at
most $b$ is the q-binomial coefficient $\qb{a+b}{a}$.  The fourth
assertion follows from the second and Lemma \ref{VS}. \end{proof}

	Next we connect $x_k$ in type $A$ with card shuffling. First
we consider the case $k=2$. Writing $x_k=\sum c_w w$ in the group
algebra, the notation $x_k^{-1}$ will denote $\sum c_w w^{-1}$.

\begin{theorem} \label{descrip} When $W$ is the symmetric group
$S_{2n}$, the element $(x_2)^{-1}$ has the following probabilistic
interpretation:

	Step 1: Choose an even number between $1$ and $2n$ with the
probability of getting $2j$ equal to $\frac{{2n \choose
2j}}{2^{2n-1}}$. From the stack of $2n$ cards, form a second pile of
size $2j$ by removing the top $j$ cards of the stack, and then putting
the bottom $j$ cards of the first stack on top of them.

	Step 2: Now one has a stack of size $2n-2j$ and a stack of
size $2j$. Drop cards repeatedly according to the rule that if stacks
$1,2$ have sizes $A,B$ at some time, then the next card comes from
stack $1$ with probability $\frac{A}{A+B}$ and from stack 2 with
probability $\frac{B}{A+B}$. (This is equivalent to choosing uniformly
at random one of the ${2n \choose 2j}$ interleavings preserving the
relative orders of the cards in each stack).

	 The description of $x_2^{-1}$ is the same for the symmetric
group $S_{2n+1}$, except that at the beginning of Step 1, the chance
of getting $2j$ is $\frac{{2n+1 \choose 2j}}{2^{2n}}$ and at the
beginning of Step 2, one has a stack of size $2n+1-2j$ and a stack of
size $2j$. \end{theorem}

\begin{proof} We argue for the case $S_{2n}$, the case of $S_{2n+1}$
being similar. Recall that in type $A_{2n-1}$ the coroot lattice is
all vectors with integer components and zero sum with respect to a
basis $e_1,\cdots,e_{2n}$, that $\alpha_i=e_i-e_{i+1}$ for
$i=1,\cdots,2n-1$ and that $\alpha_0=e_{2n}-e_1$. The elements of the
coroot lattice contributing to some $a_{2,I}$ are:

\[ \begin{array}{ll} (0,0,\cdots,0,0) & I=\tilde{\Pi}-\alpha_0\\
(1,0,0,\cdots,0,0,-1) & I=\tilde{\Pi}-\{\alpha_1,\alpha_{2n-1}\}\\
(1,1,0,0,\cdots,0,0,-1,-1) &
I=\tilde{\Pi}-\{\alpha_2,\alpha_{2n-2}\}\\ \cdots & \cdots\\
(1,1,\cdots,1,0,0,-1,\cdots,-1,-1) &
I=\tilde{\Pi}-\{\alpha_{n-1},\alpha_{n+1}\}\\
(1,1,\cdots,1,1,-1,-1,\cdots,-1,-1) & I=\tilde{\Pi}-\alpha_n
\end{array} \] One observes that the inverses of the permutations in
the above card shuffling description for a given $j$ contribute to
$u_I$ where

\[ I =             \left\{ \begin{array}{ll}
\tilde{\Pi}-\alpha_0 & \mbox{if $2j=0$}\\				
\tilde{\Pi}-\{\alpha_k,\alpha_{2n-k}\} & \mbox{if $2j=2 min(k,2n-k)$}\\
\tilde{\Pi}-\alpha_n & \mbox{if $2j=2n$}\end{array}
			\right.  \] The total number of such
permutations for a fixed value of $j$ is ${2n \choose 2j}$, the number
of interleavings of $2n-2j$ cards with $2j$ cards preserving the
relative orders in each pile. Since $\sum_{j=0}^n {2n \choose
2j}=2^{2n-1}$, and $\sum_{I \subseteq \tilde{\Pi}} a_{2,I} |U_I| =
2^{2n-1}$, the proof is complete. \end{proof}

	Note that when $n$ is prime and $k$ is a power of $n$, the
only contribution in the fourth formula comes from $d=1$. Using this
observation, the follow-up paper \cite{F4} shows that under these
conditions, the element $x_k^{-1}$ is the same as a $k$ riffle shuffle
followed by a cut at a uniform position. This observation (and Theorem
\ref{descrip}) suggest the following

{\it Problem 2:} Is there a useful ``physical'' description of the
elements $x_k$ in type $A$ for integer $k>2$, which renders some of
its algebraic properties more transparent? Such a description exists
for GSR riffle shuffles \cite{BD} and explains why a $k_1$ shuffle
followed by a $k_2$ shuffle is a $k_1k_2$ shuffle.

	Next we observe that for the identity conjugacy class in type
$A$, Theorem \ref{analogDMP} has the following consequence.

\begin{cor} \label{recip} For any positive integer $n$ and prime power
$q$, the number of ways (disregarding order and allowing repetition)
of writing $0$ mod $q-1$ as the sum of $n$ integers from the set
$0,1,\cdots,q-1$ is equal to the number of ways (disregarding order
and allowing repetition) of writing $0$ mod $n$ as the sum of $q-1$
integers from the set $0,1,\cdots,n-1$. \end{cor}

\begin{proof} Consider $k^{n-1}$ multiplied by the coefficient of the
identity in $x_q$ in type $A_{n-1}$. By part 2 of Theorem
\ref{typeAform1}, this is the number of ways of writing $0$ mod $n$ as
the sum of $q-1$ integers from the set $0,1,\cdots,n-1$. Theorem
\ref{analogDMP} states that this is the number of monic degree $n$
polynomials over $F_q$ with constant term 1 which factor into linear
terms. Working in the multiplicative group of $F_q$, this is clearly
the number of ways of writing $0$ mod $q-1$ as the sum of $n$ integers
from the set $0,1,\cdots,q-1$. \end{proof} 

	We remark that Corollary \ref{recip} holds for any positive
integers $n,q$. This can be seen from Lemma \ref{VS}. It independently
appeared in an invariant theoretic setting in \cite{EJP}.

	Next we reformulate Theorem \ref{analogDMP} in type $A$ in
terms of generating functions. This makes its number theoretic content
more visible, because one side is mod $n$ and the other side is mod
$k-1$. For its proof, Lemma \ref{likeVic} will be helpful. We use the
notation that $f_{m,k,i,d}$ is the coefficient of $z^m$ in
$(\frac{z^{kd}-1}{z^d-1})^{i/d}$.

\begin{lemma} \label{likeVic} The number of size $i$ aperiodic
necklaces on the symbols $\{0,1,\cdots,k-1\}$ with total symbol sum
$m$ is $1/i \sum_{d|i} \mu(d) f_{m,k,i,d}$ .  \end{lemma}

\begin{proof} This is an elementary Mobius inversion running along the
lines of a result in \cite{Rei1}. \end{proof}

\begin{theorem} \label{reformulation} Let $n_i(w)$ be the number of
$i$-cycles in a permutation $w$. Then Theorem 1 in type $A$ implies
the assertion (which we intentionally do not simplify) that for all
$n,k$,

\begin{eqnarray*}
& & \sum_{m=0 \ mod \ n} Coef. \ of \ q^m u^n t^k \
in \ \sum_{n=0}^{\infty} \frac{u^n} {(1-tq)\cdots(1-tq^n)} \sum_{w
\in S_n} t^{cd(w)} q^{maj(w)} \prod x_i^{n_i(w)}\\ & = &
\sum_{m=0 \ mod \ k-1} Coef. \ of \ q^m u^n t^k \ in \
\sum_{k=0}^{\infty} t^k \prod_{i=1}^{\infty} \prod_{m=1}^{\infty}
(\frac{1}{1-q^mx_iu^i})^{1/i \sum_{d|i} \mu(d) f_{m,k,i,d}}.
\end{eqnarray*} \end{theorem}

\begin{proof} The left hand side is equal to

\begin{eqnarray*}
& & \sum_{w \in S_n} \sum_{m=0 \ mod \ n} Coef. \ of
\ q^m t^{k-cd(w)} \ in \ \frac{1}{(1-tq)\cdots(1-tq^n)} q^{maj(w)}
\prod x_i^{n_i(w)}\\
& = & \sum_{w \in S_n} \sum_{m=0 \ mod \ n} Coef. \ of \ q^m \ in \ \qb{n+k-cd(w)-1}{n-1}
q^{maj(w)} \prod x_i^{n_i(w)},
\end{eqnarray*} where the last step uses Theorem 349 on page 280 of \cite{HW}. Note by part 3 of Theorem \ref{typeAform1} that this expression is precisely the cycle structure generating function under the measure $x_k$, multiplied by $k^{n-1}$.

	To complete the proof of the theorem, it must be shown that
the right hand side gives the cycle structure generating function for
degree $n$ polynomials over a field of $k$ elements with constant term
$1$. Let $\phi$ be a fixed generator of the multiplicative group of
the field $F_k$ of $k$ elements, and let $\tau_i$ be a generator of
the multiplicative group of the degree $i$ extension of $F_k$, with
the property that $\tau_i^{(k^i-1)/(k-1)}=\phi$. Recall Golomb's
correspondence \cite{Go} between degree $i$ polynomials over $F_k$ and
size $i$ aperiodic necklaces on the symbols $\{0,1,\cdots,k-1\}$. This
correspondence goes by taking any root of the polynomial, expressing
it as a power of $\tau_i$ and then writing this power base $k$ and
forming a necklace out of the coefficients of $1,k,k^2,\cdots,
k^{i-1}$. It is then easy to see that the norm of the corresponding
polynomial is $\phi$ raised to the sum of the necklace entries. The
result now follows from Lemma \ref{likeVic}. Note that there is no
$m=0$ term because the polynomial $z$ can not divide a polynomial with
constant term 1. \end{proof}

	It is perhaps interesting to compare the generating function
in Theorem \ref{reformulation} with a generating function of a similar
flavor.  For $w \in S_n$ let $d(w)=1+|\{i:w(i)>w(i+1), 1 \leq i \leq
n-1|$. Thus $d(w)$ is one more than the number of linear
descents. Gessel \cite{Ge} proved that

\begin{eqnarray*} & & \sum_{n=0}^{\infty} \frac{u^n}
{(1-t)(1-tq)\cdots(1-tq^n)} \sum_{w \in S_n} t^{d(w)} q^{maj(w)} \prod
x_i^{n_i(w)}\\ & = & \sum_{k=1}^{\infty} t^k \prod_{i=1}^{\infty}
\prod_{m=0}^{\infty} (\frac{1}{1-q^mx_iu^i})^{1/i \sum_{d|i} \mu(d)
f_{m,k,i,d}}.
\end{eqnarray*}

	This raises the following

{\it Problem 3:} What is the joint generating function for
permutations by cyclic descents, major index, and cycle structure? Can
it be used to resolve Statement 1 in Section 5 of \cite{F4}?

	We remark that Theorem 8 in \cite{F4} is equivalent to a
generating function for permutations by cyclic descents and cycle
strcuture.

\section{Type $C$ Affine Shuffles} \label{typeC}

	This section studies the $x_k$ in type $C_n$. We recall that
the elements of $C_n$ can be viewed as signed permutations $w$ on the
symbols $1,\cdots,n$. From the description of the root system of page
42 of \cite{Hu2}, it follows that (ordering the integers $1<2<3<
\cdots < \cdots < -3<-2<-1$ as in \cite{Rei1})

\begin{enumerate}

\item $w$ has a descent at position $i$ for $1 \leq i \leq n-1$ if
$w(i)>w(i+1)$.

\item $w$ has a descent at position $n$ if $w(n)<0$.

\item $w$ has a cyclic descent at position $1$ if $w(1)>0$.

\end{enumerate} For example the permutation $3 \ 1 \ -2 \ 4 \ 5$ has a
cyclic descent at position 1 and descents at positions 1 and 3.

	Lemma \ref{com}, which follows easily from Theorem 1 of
\cite{Ce2}, gives a formula for $x_k$.

\begin{lemma} \label{com} Let $d(w)$ and $cd(w)$ denote the number of
descents and cyclic descents of $w \in C_n$. Then the coefficient of
$w$ in $x_k$ is

\[ \begin{array}{ll}
\frac{1}{k^n} {\frac{k-1}{2}+n-d(w) \choose n} & \mbox{$k$ odd}\\
\frac{1}{k^n} {\frac{k}{2}+n-cd(w) \choose n} & \mbox{$k$ even}
\end{array}\] \end{lemma}

\begin{proof} For the first assertion, from Theorem 1 of \cite{Ce2},
the coefficient of $w$ in $x_k$ is \begin{eqnarray*} \frac{1}{k^n}
\sum_{l=d(w)}^n {\frac{k-1}{2} \choose l} {n-d(w) \choose l-d(w)} & =
& \frac{1}{k^n} \sum_{l=d(w)}^n {\frac{k-1}{2} \choose l} {n-d(w)
\choose n-l}\\ & = & \frac{1}{k^n} \sum_{l=0}^n {\frac{k-1}{2} \choose
l} {n-d(w) \choose n-l}\\ & = & \frac{1}{k^n} {\frac{k-1}{2}+n-d(w)
\choose n}.  \end{eqnarray*} The second assertion is similar and
involves two cases. \end{proof}

	Proposition \ref{rel} shows that the elements $x_k$ in type
$C$ arise from physical models of card-shuffling (a careful reading of
\cite{Ce2} suggests that Cellini knew this for $k=2$). The models
which follow were considered previously in the literature for the
special cases $k=2$ in \cite{BD} and for $k=3$ (and implicitly for
higher odd $k$) in \cite{BB}. The higher $k$ models and the implied
formulas for card shuffling resulting from combining Propositions
\ref{com} and \ref{rel} may be of interest (no formula is given for
the $k=2$ case in \cite{BD}). 

\begin{prop} \label{rel} The element $x_k^{-1}$ in type $C_n$ has the
following description:

	Step 1: Start with a deck of $n$ cards face down. Choose
numbers $j_1,\cdots,j_k$ multinomially with the probability of getting
$j_1,\cdots,j_k$ equal to $\frac{{n \choose
j_1,\cdots,j_k}}{k^n}$. Make $k$ stacks of cards of sizes
$j_1,\cdots,j_k$ respectively. If $k$ is odd, then flip over the even
numbered stacks. If $k$ is even, the flip over the odd numbered
stacks.

	Step 2: Drop cards from packets with probability proportional
to packet size at a given time. Equivalently, choose uniformly at
random one of the ${n \choose j_1,\cdots,j_k}$ interleavings of the
packets.
\end{prop}

\begin{proof} The proof proceeds in several cases, the goal being to
show that the inverse of the above processes generate $w$ with the
probabilities in Lemma \ref{com}. We give details for one subcase--the
others being similar--namely even $k$ when $w$ satisfies
$cd(w)=d(w)$. (The other case for $k$ even is $cd(w)=d(w)+1$). The
inverse of the probabilistic description in the theorem is as follows:

	Step 1: Start with an ordered deck of $n$ cards face
down. Successively and independently, cards are turned face up and
dealt into one of $k$ uniformly chosen random piles. The odd piles are
then flipped over (so that the cards in these piles are face down).

	Step 2: Collect the piles from pile 1 to pile $k$, so that
pile 1 is on top and pile $k$ is on the bottom.

Consider for instance the permutation $w$ given in 2-line form by $-2\
3\ 1\ 4\ -6\ -5\ 7$. Note that this satisfies $cd(w)=d(w)$ because the
top card has a negative value (i.e. is turned face up). It is
necessary to count the number of ways that $w$ could have arisen from
the inverse description. This one does using a bar and stars argument
as in \cite{BD}. Here the stars represent the $n$ cards, and the bars
represent the $k-1$ breaks between the different piles. It is easy to
see that each descent in $w$ forces the position of two bars, except
for the first descent which only forces one bar. Then the remaining
$(k-1)-(2d(w)-1)=k-2d(w)$ bars must be placed among the $n$ cards as
$\frac{k-2d(w)}{2}$ consecutive pairs (since the piles alternate
face-up, face-down). This can be done in ${\frac{k}{2}+n-cd(w) \choose
n}$ ways, proving the result. \end{proof}

	We remark that Proposition \ref{rel} leads to a direct proof
of the convolution property in type $C$.

	Next recall the notion of total variation distance
$||P_1-P_2||$ between two probability distributions $P_1$ and $P_2$ on
a finite set $X$. It is defined as \[ \frac{1}{2} \sum_{x \in X}
|P_1(x)-P_2(x)|.\] The book \cite{D} explains why this is a natural
and useful notion of distance between probability distributions. The
remainder of this section computes the total variation distance of an
Affine type $C$ $k$-shuffle to uniform in the case that $k$ is
even. The paper \cite{BD} attributes an equivalent result to Bob Beals
when $k$ is a power of 2 (unpublished), but with a quite different
method of proof. We omit the case of odd $k$ as the convergence rate
to randomness has been determined in \cite{BB}.

\begin{lemma} Let $N_r$ be the number of $w$ in $C_n$ with $r$ cyclic
descents. Let $A_r$ be the number $w$ in $S_n$ with $r$
descents. Then $N_{r+1}=2^n A_r$. \end{lemma}

\begin{proof} Lemma \ref{com} shows that the chance that an Affine
Type $C$ $k$ shuffle gives a signed permutation $w$ is

\[ \frac{1}{k^n} {k/2 + n - cd(w) \choose n}.\] Using the fact that
these shuffles are a probability measure and dividing both sides of
the resulting equation by $2^n$, it follows that

\[ \sum_{r=1}^n \frac{N_r}{2^n} {k/2+n-r \choose n} = (k/2)^n.\] This
can be rewritten as

\[ \sum_{r=0}^{n-1} \frac{N_{r+1}}{2^n} {k/2+n-r-1 \choose n} =
(k/2)^n.\] Since this is true for all $k$, the relation can be
inverted to solve for $N_{r+1}$. In the theory riffle shuffles
\cite{BD} one gets the equation (Worpitzky's identity)

\[ \sum_{r=0}^{n-1} A_r {k+n-r-1 \choose n} = k^n\] for all
$k$. Thus $N_{r+1} = 2^n A_r$ as desired. \end{proof}

\begin{theorem} \label{TV} The total variation distance of an Affine
type $C$ $k$ shuffle with $k$ even to uniform is equal to the total
distance of a Bayer-Diaconis $k/2$ riffle shuffle on $S_n$ to
uniform. \end{theorem}

\begin{proof} Lemma \ref{com} shows that the chance that an Affine
type $C$ $k$ shuffle gives a signed permutation $w$ is

\[ \frac{1}{k^n} {k/2 + n - cd(w) \choose n}.\] Thus the total
variation distance is equal to

\begin{eqnarray*}
& & \sum_{r=0}^{n-1} N_{r+1} \left|\frac{1}{k^n} {k/2+n-r-1 \choose n} -
\frac{1}{2^nn!}\right|\\
& = & \sum_{r=0}^{n-1} 2^n A_r \left|\frac{1}{k^n} {k/2+n-r-1 \choose n} -
\frac{1}{2^nn!}\right|\\
& = & \sum_{r=0}^{n-1} A_r \left|\frac{1}{(k/2)^n} {k/2+n-r-1 \choose n} -
\frac{1}{n!}\right|.
\end{eqnarray*} From \cite{BD}, one recognizes this last expression as the
total variation distance between a $k/2$ riffle shuffle and
uniform. \end{proof}

\section{Dynamical Systems} \label{dynamical}

	Much of this section relates to the enumeration of unimodal
permutations by cycle structure. This problem is given two solutions:
one using a more fundamental result of Gannon \cite{Ga} and symmetric
functions, and another using Theorem \ref{analogDMP} and the
factorization theory of palindromic polynomials (which actually proves
a more general result). Some asymptotic consequences are derived. We
give a new proof of a result of Reiner \cite{Rei1}.

	A unimodal permutation $w$ on the symbols $\{1,\cdots,n\}$ is
defined by requiring that there is some $i$ with $1 \leq i \leq n$
such that the following two properties hold:
 \begin{enumerate} \item If $a<b\leq i$, then $w(a)<w(b)$.
  \item If $i \leq a<b$, then $w(a)>w(b)$.
  \end{enumerate} Thus $i$ is where the maximum is achieved, and the
permutations $12\cdots n$ and $nn-1\cdots 1$ are counted as
unimodal. For each fixed $i$ there are ${n-1 \choose i-1}$ unimodal
permutations with maximum $i$, hence a total of $2^{n-1}$ such
permutations. As noted in \cite{Ga}, unimodal permutations are those
which avoid the patterns $213$ and $312$.

	Motivated by biology and dynamical systems, Rogers \cite{Ro}
posed the problem of counting unimodal permutations by cycle
structure. This problem was solved by Gannon who gave a constructive
proof of the following elegant (and more fundamental) result. For its
statement, one defines the shape $s$ of a cycle $(i_1 \cdots i_k)$ on
some $k$ distinct symbols (call them $\{1,\cdots,k\}$) to be the cycle
$(\tau(i_1) \cdots \tau(i_k))$ where $\tau$ is the unique order
preserving bijection between $\{i_1,\cdots,i_k\}$ and
$\{1,\cdots,k\}$. Thus the shape of $(523)$ is $(312)$.

\begin{theorem} \label{Gannon} (\cite{Ga}) Let $s_1,s_2,\cdots$ denote
the possible shapes of transitive unimodal permutations. Then the
number of unimodal permutations with $n_i$ cycles of shape $s_i$ is
$2^{l-1}$, where $l$ is the number of $i$ for which $n_i>0$. \end{theorem}

	Theorem \ref{Gannon} can be rewritten in terms of generating
functions.

\begin{cor} \label{cindex} Let $n_s(w)$ be the number of cycles of $w$
of shape $s$. Let $|s|$ be the number of elements in $s$. Then

\[ 1+\sum_{n=1}^{\infty} \frac{u^n}{2^{n-1}} \sum_{w \in S_n \atop w \
unimodal} \prod_{s \ shape} x_s^{n_s(w)} = \prod_{s \ shape}
(\frac{2^{|s|}+x_su^{|s|}}{2^{|s|}-x_su^{|s|}})\]

\[ (1-u)+\sum_{n=1}^{\infty} \frac{(1-u) u^n}{2^{n-1}} \sum_{w \in S_n
\atop w \ unimodal} \prod_{s \ shape} x_s^{n_s(w)} = \prod_{s \ shape}
(\frac{2^{|s|}+x_su^{|s|}}{2^{|s|}+u^{|s|}})(\frac{2^{|s|}-u^{|s|}}{2^{|s|}
-x_su^{|s|}})\] \end{cor}

\begin{proof} For the first equation, consider the coefficient of
$\prod_s x_s^{n_s} u^{\sum |s| n_s}$ on the left hand side. It is the
probability that a uniformly chosen unimodal permutation on $\sum |s|
n_s$ symbols has $n_s$ cycles of shape $s$. The coefficient on the
right hand side is $2^{|\{s:n_s>0\}|-n}$. These are equal by Theorem
\ref{Gannon}. To deduce the second equation, observe that setting all
$x_s=1$ in the first equation gives that \[ \frac{1}{1-u} = \prod_{s \
shape} \frac{2^{|s|}+u^{|s|}}{2^{|s|}-u^{|s|}} .\] Taking reciprocals
and multiplying by the first equation yields the second
equation. \end{proof}

	The second equation in Corollary \ref{cindex} has an
attractive probabilistic interpretation. Fix $u$ such that
$0<u<1$. Then choose a random symmetric group so that the chance of
getting $S_n$ is equal to $(1-u)u^n$. Choose a unimodal $w \in S_n$
uniformly at random. Then the random variables $n_s(w)$ are
independent, each having distribution a convolution of a
binomial$(\frac{u^{|s|}} {2^{|s|}+u^{|s|}})$ with a
geometric$(1-\frac{u^{|s|}}{2^{|s|}})$.

	As another illustration of the second equation in Corollary
\ref{cindex}, we deduce the following corollary, extending the
asymptotic results in \cite{Ga} that asymptotically $2/3$ of all
unimodal permutations have fixed points and $2/5$ have $2$-cycles.

\begin{cor} \label{bigeom} In the $n \rightarrow \infty$ limit, the
random variables $n_s$ converge to the convolution of a
binomial$(\frac{1} {2^{|s|}+1})$ with a
geometric$(1-\frac{1}{2^{|s|}})$ and are asymptotically independent.
\end{cor}

\begin{proof} The result follows from the claim that if $f(u)$ has a
Taylor series around 0 and $f(1)<\infty$, then the $n \rightarrow
\infty$ limit of the coefficient of $u^n$ in $\frac{f(u)}{1-u}$ is
$f(1)$. To verify the claim, write the Taylor expansion $f(u) =
\sum_{n=0}^{\infty} a_n u^n$ and observe that the coefficient of $u^n$
in $\frac{f(u)}{1-u} = \sum_{i=0}^n a_i$.  \end{proof}

	Rogers and Weiss \cite{RW} used dynamical systems to count the
number of transitive unimodal permutation on $n$ symbols. We offer a
proof using symmetric function theory. The paper \cite{T} generalizes
this argument to arbitrary cycle types.

	Some notation is needed. A subset $D=\{d_1,\cdots,d_k\}$ of
$\{1,2,\cdots,n-1\}$ defines a composition $C(D)$ of $n$ with parts
$d_1,d_2-d_1,\cdots,n-d_k$. A standard Young tableau is said to have a
descent at position $i$ if $i+1$ occurs in a row lower than $i$. The
descent set of a standard Young tableau thus defines a composition of
$n$.

\begin{lemma} \label{countem} The number of transitive unimodal
permutations on $n$ symbols is 

\[\frac{1}{2n} \sum_{d|n \atop d \ odd} \mu(d) 2^{\frac{n}{d}}.\]
\end{lemma}

\begin{proof} Symmetric function notation from Chapter 1 of Macdonald
\cite{Mac} is used. Thus $p_{\lambda},h_{\lambda}$ $,
e_{\lambda},s_{\lambda}$ are the power sum, complete, elementary, and
Schur symmetric functions parameterized by a partition $\lambda$. From
Theorem 2.1 of \cite{G}, the number of $n$ cycles with descent set $D$
is the inner product of a Lie character $L_n=\frac{1}{n} \sum_{d|n}
\mu(d) p_d^{\frac{n}{d}}$ and a Foulkes character $F_{C(D)}$. From the
proof of Corollary 2.4 of \cite{G}, $F_{C(D)}=\sum_{|\lambda|=n}
\beta_{\lambda} s_{\lambda}$ where $\beta_{\lambda}$ is the number of
standard tableaux of shape $\lambda$ with descent composition
$C(D)$. Thus the sought number is \[ <\frac{1}{n} \ \sum_{d|n} \mu(d)
p_d^{\frac{n}{d}}, e_n+\sum_{i=2}^{n-1} s_{i,(1)^{n-i}}+h_n>.\]
Expanding these Schur functions using exercise 9 on page 47 of
\cite{Mac}, using the fact that the $p_{\lambda}$ are an orthogonal
basis of the ring of symmetric functions with known normalizing
constants (page 64 of \cite{Mac}), and using the expansions of $e_n$
and $h_n$ in terms of the $p_{\lambda}$'s (page 25 of \cite{Mac}) it
follows that

\begin{eqnarray*}
 & & <\frac{1}{n} \ \sum_{d|n} \mu(d) p_d^{\frac{n}{d}},
e_n+\sum_{i=2}^{n-1} s_{i,(1)^{n-i}}+h_n>\\
&=& <\frac{1}{n} \ \sum_{d|n} \mu(d) p_d^{\frac{n}{d}},
\sum_{i \ even} h_i e_{n-i}>\\
&=& \frac{1}{n} \sum_{d|n} \mu(d) <p_d^{\frac{n}{d}},\sum_{i=1,\cdots,\frac{n}{d} \atop di \ even} h_{di} e_{n-di}>\\
&=& \frac{1}{n} \sum_{d|n} \mu(d) <p_d^{\frac{n}{d}}, p_d^{\frac{n}{d}} \sum_{i=1,\cdots,\frac{n}{d} \atop di \ even} \frac{(-1)^{n-di-\frac{n}{d}+i}}{d^{\frac{n}{d}} i! (\frac{n}{d}-i)!}>\\
& = & \frac{1}{n} \sum_{d|n} \mu(d) (-1)^{n-\frac{n}{d}} \sum_{i=1,\cdots,\frac{n}{d} \atop di \ even} (-1)^i {\frac{n}{d} \choose i}\\
& = & \frac{1}{2n} \sum_{d|n \atop d \ odd} \mu(d) 2^{\frac{n}{d}}.
\end{eqnarray*}
\end{proof}

	Corollary \ref{cindex} and Lemma \ref{countem} have the
following immediate consequence.

\begin{cor} Let $n_i(w)$ be the number of $i$-cycles of a permutation
$w$. Then \[ 1+\sum_{n=1}^{\infty} \frac{u^n}{2^{n-1}} \sum_{w \in S_n
\atop w \ unimodal} \prod_i x_i^{n_i(w)} = \prod_i
(\frac{2^i+x_iu^i}{2^i-x_iu^i})^{\frac{1}{2i} \sum_{d|i \atop d \ odd}
\mu(d) 2^{\frac{i}{d}}}.\] \end{cor}

	Theorem \ref{analogDMP} will yield a second proof of the
enumeration of unimodal permutations by cycle structure by relating
the problem to the factorization theory of palindromic polynomials
over finite fields. The first step is to reformulate Theorem
\ref{analogDMP} in type $C$. The following lemmas, the first of which
is well known, will be helpful. The symbol $\mu$ denotes the Moebius
function of elementary number theory.

\begin{lemma} \label{totalnum} The number of degree $n$ irreducible
polynomials over $F_q$ is equal to \[ \frac{1}{n} \sum_{d|n} \mu(d)
q^{n/d}.\] \end{lemma}

\begin{lemma} \label{countinvar} (\cite{FNP}) Let $e=1$ if $q$ is even
and $e=2$ if $q$ is odd. Then the number of monic, degree $n$
polynomials $f(z)$ over $F_q$ with non-zero constant coefficient and
invariant under the involution $f(z) \mapsto f(0)^{-1} z^n
f(\frac{1}{z})$ is

\[ \left\{ \begin{array}{ll}
																					e & \mbox{if
$n=1$}\\
																																			0 & \mbox{if $n$ \ is \ odd \ and \ $n>1$}\\
																						\frac{1}{n} \sum_{d|n \atop d \ odd} \mu(d) (q^{\frac{n}{2d}}+1-e)		& \mbox{Otherwise}
																																				\end{array}
			\right.			 \]

\end{lemma}

	 Recall that $x_q(w)$ denotes the coefficient of $w$ in $x_q$.

\begin{theorem} \label{reduce} Let $e=1$ if $q$ is even and $e=2$ if
$q$ is odd. Let $\lambda_i(w)$ and $\mu_i(w)$ be the number of
positive and negative $i$-cycles of a signed permutation $w$ in
$C_n$. Then

\begin{eqnarray*}
& & 1+\sum_{n \geq 1} u^nq^n \sum_{w \in C_n} x_q(w) \prod_{i \geq 1}
x_i^{\lambda_i(w)} y_i^{\mu_i(w)}\\
& = & (\frac{1}{1-x_1u})^{e-1} \prod_{m
\geq 1} (\frac{1+y_m u^m}{1-x_m u^m})^{\frac{1}{2m} \sum_{d|m \atop d
\ odd} \mu(d) (q^{\frac{m}{d}}+1-e)}.
\end{eqnarray*}

\end{theorem}

\begin{proof} One argues separately for odd and even characteristic
and first for prime powers. Taking the coefficient of $u^n \prod_i
x_i^{\lambda_i} y_i^{\mu_i}$ on the left hand side of this equation
and dividing by $q^n$ gives by Lemma \ref{com} the probability that
$w$ chosen according to the $x_q$ probability measure is in a
conjugacy class with $\lambda_i$ positive $i$-cycles and $\mu_i$
negative $i$-cycles for each $i$. By Theorem \ref{analogDMP}, to
verify the theorem for even prime powers it is enough to show that the
coefficient of $u^n \prod_i x_i^{\lambda_i} y_i^{\mu_i}$ on the right
hand side of this equation is the number of degree $2n$ monic
self-conjugate polynomials over $F_q$ which factor as

\[ \prod_{\{\phi_j, \bar{\phi_j}\}} [\phi_j \bar{\phi_j}]^{r_{\phi_j}}
\prod_{\phi_j : \phi_j = \bar{\phi_j}} \phi_j^{s_{\phi_j}} \] (with
$\phi_j$ where $s_{\phi_j} \in \{0,1\}$) and $\lambda_i=\sum_{\phi:
deg(\phi)=i} r_{\phi}$ and $\mu_i=\sum_{\phi: deg(\phi)=2i}
s_{\phi}$. This follows readily from Lemmas \ref{totalnum} and
\ref{countinvar}. The theorem now follows for arbitrary $q$ since two
functions analytic in a region and agreeing on a set with an
accumulation point ($q=\infty$) in that region must be
equal. \end{proof}

	Corollary \ref{re} deduces the enumeration of unimodal
permutations by cycle structure.

\begin{cor} \label{re} Let $n_i(w)$ be the number of $i$-cycles of
$w \in S_n$. Then

\[ 1+\sum_{n=1}^{\infty} \frac{u^n}{2^{n-1}} \sum_{w \in S_n \atop w \
unimodal} \prod_i x_i^{n_i(w)} = \prod_i
(\frac{2^i+x_iu^i}{2^i-x_iu^i})^{\frac{1}{2i} \sum_{d|i \atop d \ odd}
\mu(d) 2^{\frac{i}{d}}}.\] \end{cor}

\begin{proof} Given Theorem \ref{reduce} with $q=2$, it is enough to
define a $2$ to $1$ map $\eta$ from the $2^n$ type $C_n$
characteristic 2 shuffles to unimodal elements of $S_n$, such that
$\eta$ preserves the number of $i$-cycles for each $i$, disregarding
signs. To define $\eta$, recalling Proposition \ref{rel} observe that
the $2$ shuffles are all ways of cutting a deck of size $n$, then
flipping the first pack, and choosing a random interleaving. For
instance if one cuts a 12 card deck at position 6, such an
interleaving could be \[[-6,-5,7,8,-4,9,-3,10,-2,11,-1,12].\] Observe
that taking the inverse of this permutation and disregarding signs
gives \[[11,9,7,5,2,1,3,4,6,8,10,12].\] Next one conjugates by the
involution transposing each $i$ with $n+1-i$, thereby obtaining a
unimodal permutation. Note that this map preserves cycle structure,
and is $2$ to $1$ because the first symbol (in the example $-6$, can
always have its sign reversed yielding a possible shuffle).
\end{proof}

	The following corollary describes the $n \rightarrow \infty$
asymptotics of cycle structure for type $C$ affine $q$ shuffles. We
omit the proof, which is essentially the same as Corollary
\ref{bigeom}.

\begin{cor} \label{asymptot} Let $\lambda_i(w)$ and $\mu_i(w)$ be the
number of positive and negative $i$-cycles of a signed permutation $w$
in $C_n$.

\begin{enumerate}

\item Fix $u$ such that $0<u<1$. Then choose a random hyperoctahedral
group so that the chance of getting $C_n$ is equal to
$(1-u)u^n$. Choose $w \in C_n$ according to the affine $q$ shuffle
measure. Then the random variables $\{\lambda_m,\mu_m\}$ are
independent. The $\lambda_m$ ($m \geq 2$) are distributed as the
convolution of $\frac{1}{2m} \sum_{d|m} \mu(d) (q^{m/d}+1-e)$ many
geometrics with parameter $1-\frac{u^m}{q^m}$ and $\lambda_1$ is
distributed as the convolution of $\frac{1}{2} (q+e-1)$ many
geometrics with parameter $1-\frac{u}{q}$. The $\mu_m$ are distributed
as the convolution of $\frac{1}{2m} \sum_{d|m} \mu(d) (q^{m/d}+1-e)$
many binomials with parameter $\frac{u^m/q^m}{1+u^m/q^m}$.

\item Choose $w \in C_n$ according to the affine $q$ shuffle
measure. Then in the $n \rightarrow \infty$ limit, the random
variables $\{\lambda_m,\mu_m\}$ are independent. The $\lambda_m$ ($m
\geq 2$) are distributed as the convolution of $\frac{1}{2m}
\sum_{d|m} \mu(d) (q^{m/d}+1-e)$ many geometrics with parameter
$1-\frac{1}{q^m}$ and $\lambda_1$ is distributed as the convolution
of $\frac{1}{2} (q+e-1)$ many geometrics with parameter
$1-\frac{1}{q}$. The $\mu_m$ are distributed as the convolution of
$\frac{1}{2m} \sum_{d|m} \mu(d) (q^{m/d}+1-e)$ many binomials with
parameter $\frac{1/q^m}{1+1/q^m}$.

\end{enumerate}
\end{cor}

{\bf Remark:} Type $C_n$ shuffles also relate to dynamical systems in
another way, analogous to the type $A$ construction for Bayer-Diaconis
shuffles \cite{BD}. Here we describe
the case $k=2$. One drops $n$ points in the interval $[-1,1]$
uniformly and indepedently. Then one applies the map $x \mapsto
2|x|-1$. The resulting permutation can be thought of as a signed
permutation, since some points preserve and some reverse
orientation. From Proposition \ref{rel}, this signed permutation
obtained after iterating this map $r$ times has the distribution of
the type $C_n$ shuffle with $k=2^r$. Lalley \cite{La1} studied the
cycle structure of random permutations obtained by tracking $n$
uniformly dropped points after iterating a map a large number of
times. His results applied to piecewise monotone maps, and he proved
that the limiting cycle structure is a convolution of
geometrics. Hence Corollary \ref{asymptot} shows that Lalley's results
do not extend to functions such as $x \mapsto 2|x|-1$.

	As a final result, we deduce a new proof of the following
result of Reiner. Here $d(w)$ denotes the number of descents of $w \in
C_n$.

\begin{cor} (\cite{Rei1})

\begin{eqnarray*}
& & \sum_{n \geq 0} \frac{u^n}{(1-t)^{n+1}} \sum_{w \in C_n} t^{d(w)+1}
\prod_i x_i^{\lambda_i(w)} y_i^{\mu_i(w)}\\
& = & \sum_{k \geq 0} t^k
\frac{1}{1-x_1u} \prod_{m \geq 1}
(\frac{1+y_mu^m}{1-x_mu^m})^{\frac{1}{2m} \sum_{d|m \atop m \ odd}
\mu(d) ((2k-1)^{m/d}-1)}.
\end{eqnarray*}
\end{cor}

\begin{proof} Taking coefficients of $t^k$ on both sides of the
equation under question and then setting $q=2k-1$ gives the equation

\[ \sum_{n \geq 0} u^n \sum_{w \in C_n} {\frac{q-1}{2}+n-d(w) \choose
n} \prod_i x_i^{\lambda_i(w)} y_i^{\mu_i(w)} = \frac{1}{1-x_1u}
\prod_{m \geq 1} (\frac{1+y_mu^m}{1-x_mu^m})^{\frac{1}{2m} \sum_{d|m
\atop m \ odd} \mu(d) (q^{m/d}-1)}.\] However this equation follows
from Theorem \ref{reduce} for odd $q$ and Lemma \ref{com}. \end{proof}

\section{Proof of Theorem \ref{analogDMP}: The Brauer Complex}
\label{refine}

	The purpose of this section is to report a proof, due to
Professor Roger Carter, of Theorem \ref{analogDMP}. The proof uses a
geometric object called the Brauer complex. All relevant background
(including pictures) can be found in Section 3.8 of \cite{C}.
	
	Let $Y$ be the coroot lattice and $W$ the Weyl group, so that
$<Y,W>$ is the affine Weyl group. The group $<Y,W>$ acts on the vector
space $Y \otimes R$ with $Y$ acting by translations $T_y:v \mapsto
v+y$ and $W$ acting by orthogonal transformations. The affine Weyl
group has a fundamental region in $Y \otimes R$ given by \[ \bar{A_1}
= \{v \in Y \otimes R| <\alpha_i,v> \geq 0 \ for \ i=1,\cdots,r,
<-\alpha_0,v> \leq 1\}.\] Let $Q_{p'}$ be the additive group of
rational numbers $\frac{s}{t}$ where $s,t \in Z$ and $t$ is not
divisible by $p$ (the characteristic). Proposition 3.8.1 of \cite{C}
shows that there is an action of $F$ on $\bar{A_{p'}}=\bar{A_1} \cap
(Y \otimes Q_{p'})$ given by taking the image of $v \in \bar{A_{p'}}$
to be the unique element of $\bar{A_{p'}}$ equivalent to $F(v)$ under
$<Y,W>$.

	We highlight the following facts.

{\it Fact 1:} (Proposition 3.7.3 of \cite{C}) The $q^r$ semisimple
conjugacy classes of $G^F$ are in bijection with the $q^r$ elements of
$\bar{A_{p'}}$ which are stable under the action of $F$.

	As an example, for type $A_2$ with $q=3$, the 9 stable points
are

\begin{enumerate}

\item $(0,0,0)$,$(1/2,0,-1/2)$ corresponding to polynomials that
factor into linear pieces

\item $(1/4,0,-1/4),(4/8,1/8,-5/8),(5/8,-1/8,-4/8)$ corresponding
to polynomials that are a product of a linear and a degree 2 factor

\item
$(6/26,2/26,-8/26),(8/26,-2/26,-6/26),(10/26,4/26,-14/26),(14/26,-4/26,-10/26)$
corresponding to irreducible polyomials.

\end{enumerate} The point $v=(4/8,1/8,-5/8)$ for instance is stable
because $3v=v^{(23)}+(1,1,-2)$.

	We remark that this bijection is not entirely canonical,
because the isomorphism between the multiplicative group of the
algebraic closure of $F_q$ and $Q_{p'}/Z$ (Proposition 3.1.3 of
\cite{C}) is not entirely canonical. In other words, we have chosen (a
consistent set of) generators of the multiplicative groups of all of
the finite extensions of $F_q$. The bijection in Fact 1 is canonical
only after this choice.

{\it Fact 2:} (Corollary 3.8.3 of \cite{C}) There is a bijection
between semisimple conjugacy classes in $G^F$ and simplices of maximal
dimension in the Brauer complex.

	For all that follows,

	\[ \bar{A_q} = \{ v \in Y \otimes R: <\alpha_i,v> \geq 0 \ for \
i=1,\cdots,r, <-\alpha_0,v> \leq q\} \] Let $I(y)$ be the set of
$\alpha_i$ with $i \in \{0,1,\cdots,r\}$ such that $y$ lies on the
$i$-boundary wall of $\bar{A_q}$.

	We now prove Theorem \ref{analogDMP}.

\begin{proof} The proof proceeds in two steps. Step 1 is to show that
there is a bijection between semisimple classes $c$ in $G^F$ and pairs
$(y,w) \in Y \times W$ such that $y \in Y \cap \bar{A_q}$ and $I(y) \cap
Cdes(w^{-1}) = \emptyset$. Step 2 is to show that $\Phi(c)$ is
conjugate to $w$. The theorem then follows from the definition of
$x_q$.

	{\it Step 1:} It is known that there is a bijection between
simplices of maximal dimension in the Brauer complex and elements
$\omega$ in the affine Weyl group such that $\omega(\bar{A_1}) \subset
\bar{A_q}$. The alcoves $\omega(\bar{A_1})$ are all obtained by first
transforming by $w \in W$ to give $w(\bar{A_1})$ and then translating
by $T_y$ for some $y \in Y$. Let $S$ be the union of the alcoves
$w(\bar{A_1})$ for $w \in W$; $S$ is called the basic star. The sets
$T_y(S)$ are called the stars and the centers of the stars are the
elements of $Y$.

	Each alcove $\omega(\bar{A_1})$ which lies in $\bar{A_q}$ lies
in some star whose center lies in $Y \cap \bar{A_q}$. Conversely, if
$y \in Y \cap \bar{A_q}$ we wish to know which alcoves in the star
$T_y(S)$ lie in $\bar{A_q}$. If $y$ does not lie on any boundary wall
of $\bar{A_q}$ all alcoves in $T_y(S)$ lie in $\bar{A_q}$. If $y$ lies
on the boundary wall corresponding to $i \in \{0,1,\cdots,r\}$ then
the alcove $T_yw(\bar{A_1})$ lies on the $\bar{A_q}$-side of this
boundary wall if and only if $w^{-1}(\alpha_i)$ is a positive
root. This can be seen by looking at the star $S$. Thus there is a
bijection between semisimple classes $c$ in $G^F$ and pairs $(y,w) \in
Y \times W$ such that $y \in Y \cap \bar{A_q}$ and $I(y) \cap
Cdes(w^{-1}) = \emptyset$.

	{\it Step 2:} Let $T_0$ be a maximal split torus of $G$ and
let $Y_0=Hom(Q_{p'}/Z,T_0)$ be its co-character group. Let $T$ be an
$F$-stable maximal torus of $G$ obtained from $T_0$ by twisting with
$w \in W$. We have conjugation maps $T \mapsto T_0, Y \mapsto
Y_0$. Under these maps $F:Y \mapsto Y$ maps to $w^{-1}F: Y_0 \mapsto
Y_0$. 

	Let $\omega$ be an element of the affine Weyl group such that
$\omega(\bar{A_1}) \subset \bar{A_q}$. From Section 3.8 of \cite{C},
$\bar{A}$ contains a unique $p$ satisfying $F^{-1}\omega(p)=p$
i.e. $F(p)=p^w+y_0$. Let the walls of $\bar{A_q}$ be
$H_0,H_1,\cdots,H_n$. Let $J=\{i: p \in H_i\}$. $J$ is a proper subset
of $\{0,1,\cdots,n\}$. The roots $\alpha_i, i \in J$ form a simple
system $\Pi_J$ in a subsystem $\Phi_J \subset \Phi$. From page 102 of
\cite{C} the point $p$ maps to an element $t_0 \in T_0$. Then
$F(t_0)=t_0^w$ and $t_0$ lies in the semisimple conjugacy class of
$G^F$ corresponding to the point $p$. $\Phi_J$ can be identitied with
the root system of the centralizer of the semisimple conjugacy class
of $G^F$ corresponding to the point $t_0$.

	To complete the proof of Step 2, it suffices to show (by Lemma
\ref{explicit}) that $w(\Phi_J^+)=\Phi_J^+$. The construction of the
point $p$ as the intersection of a sequence of increasingly small
alcoves, each obtained from the previous one by a map $F^{-1} \omega$
which preserves the type of the walls, shows that $p$ lies in the
$J$-face of $\bar{A_q}$ and of $\omega(\bar{A_1})$ (the $J$-face of
$\bar{A_q}$ is the intersection of the $H_i$ for $i \in J$). For $i
\in J$ let the wall $H_i$ of $\bar{A_q}$ coincide with the wall of
type $j$ of $\omega(\bar{A_1})$. The root orthogonal to $H_i$ pointing
into $\bar{A_q}$ is $\alpha_i$. Consider the root orthogonal to the
wall of type $j$ for $\omega(\bar{A_1})$ pointing into
$\omega(\bar{A_1})$. Since $\omega(\bar{A_1}) =T_{y_0}(w(\bar{A_1}))$ this
is the root orthogonal to the wall of type $j$ for $w(\bar{A_1})$
pointing into $w(\bar{A_1})$. This is the root $w(\alpha_j)$ since the
wall of type $j$ for $w(\bar{A_1})$ is the image under $w$ of the wall
of type $j$ for $\bar{A_1}$. Hence $\alpha_i=w(\alpha_j)$. This shows,
since $i,j \in J$, that $w(\Pi_J)=\Pi_J$. Hence
$w(\Phi_J^+)=\Phi_J^+$, as desired.  \end{proof}

	The results of this section raise

{\it Problem 4:} Is there an analog of the Brauer complex that sheds
light on the problems in \cite{F3}, or gives a general type
Gessel/Reutenauer bijection?

\section{Acknowledgements} This research was supported by an NSF
Postdoctoral Fellowship. The author is very indebted to Roger Carter
for proving conjectures stated in an early version of this paper,
superseding the author's partial results on the problem. Persi
Diaconis offered expository feedback and motivated the author to seek
a new proof of Beals' work on total variation distance.

\end{document}